\documentclass[12pt,reqno]{amsart}

\usepackage{setspace,hyperref}

\author[V. Kanovei]{Vladimir Kanovei} \address{V. Kanovei, IPPI,
Moscow, and MIIT, Moscow, Russia}\email{kanovei@googlemail.com}

\author[K. Katz]{Karin U. Katz}\address{K. Katz, Department of
Mathematics, Bar Ilan University, Ramat Gan 52900 Israel}
\email{katzmik@math.biu.ac.il}

\author[M. Katz]{Mikhail G. Katz}\address{M. Katz, Department of
Mathematics, Bar Ilan University, Ramat Gan 52900 Israel}
\email{katzmik@macs.biu.ac.il}

\author[D. Sherry]{David Sherry}\address{D. Sherry, Department of
Philosophy, Northern Arizona University, Flagstaff, AZ 86011, US}
\email{David.Sherry@nau.edu}

\begin{document}


\thispagestyle{empty}

\title{Euler's lute and Edwards' oud}

\keywords{Infinitesimals; ratio of vanishing increments \quad\quad\quad\quad
{\large Published in $\phantom{I}^{\phantom{I}}$
\emph{The Mathematical Intelligencer}, see \quad\quad\quad\quad
\\
\url{http://dx.doi.org/10.1007/s00283-015-9565-6}
}.}

\begin{abstract}
In a piece published in 1981, H. M. Edwards touts the benefits of
\emph{reading the masters}.  A quarter-century later, Edwards takes
seriously his own advice by publishing an encomium on Euler's
\emph{Institutiones} (1755).  While we agree with Edwards that we
shall all do well by studying the masters, we argue that to derive the
full benefit, one must read historically important texts without
interposing a lens formed by one's own mathematical accomplishments.
We show, in particular, that Edwards misses much of the style and
substance of Euler's \emph{Institutiones} by reading the text through
a constructivist lens.

\end{abstract}

\maketitle


\section{Euler on ratios of vanishing increments}
\label{one}

\emph{Read the masters!}  is Harold Edwards' text shortlisting Euler
among the \emph{masters} \cite[p.~105]{Ed81}.  Let us therefore read
Euler.  In the Preface to his 1755 \emph{Institutiones}, Euler writes:
\begin{quote}
In this way, we are led to a definition of \emph{differential
calculus}: It is \emph{a method for determining the ratio of the
vanishing increments that any functions take on when the variable, of
which they are functions, is given a vanishing
increment}. \cite[p.~vii]{Eu55b} (as translated by Blanton)
\end{quote}
Euler goes on to write that the vanishing increments involved in that
ratio
\begin{quote}
are called \emph{differentials}, and since they are without quantity,
they are also said to be \emph{infinitely small}.  Hence, by their
nature they are to be so interpreted as absolutely nothing, or they
are considered to be equal to nothing.  (ibid.)
\end{quote}
Euler's reference to \emph{equality} in his comment about
infinitesimals being ``equal to nothing'' can be interpreted in terms
of a generalized notion of equality as follows.  Euler discusses two
modes of comparison, arithmetic and geometric, in \S\,84-87 of
\cite{Eu55}.  Euler illustrates the arithmetic mode by~$ndx=0$ and the
geometric by~$\frac{a\pm ndx}{a}=1$.  Euler's generalized notions of
equality are similar to Leibniz's.  Leibniz used a generalized notion
of \emph{equality up to} in the context of his \emph{transcendental
law of homogeneity} \cite{KS2}, showing Leibnizian calculus to have
been more soundly founded than George Berkeley's criticism thereof.
Euler goes on to give an example in terms of an
infinitesimal~$\omega$:
\begin{quote}
Thus, if the quantity~$x$ is given an increment~$\omega$, so that it
becomes~$x + \omega$, its square~$x^2$ becomes~$x^2 + 2x\omega +
\omega^2$, and it takes the increment~$2x\omega+\omega^2$.  Hence, the
increment of~$x$ itself, which is~$\omega$, has the ratio to the
increment of the square, which is~$2x\omega+\omega^2$, as~$1$
to~$2x+\omega$. This ratio reduces to 1 to~$2x$, at least
when~$\omega$ vanishes. Let~$\omega = 0$, and the ratio of these
vanishing increments, which is the main concern of differential
calculus, is as 1 to~$2x$. \cite[p.~vii]{Eu55b}
\end{quote}
Euler's procedure involves an increment~$\omega$ described as being
\emph{infinitely small}.  In the first few chapters of
\emph{Institutiones}, the ratios of vanishing increments for a few
elementary functions are determined via power series expansions
obtained previously in the \emph{Introductio} \cite{Eu48}.  However,
in Chapter 4 a different picture begins to emerge, including formulas
for relations among more complicated differentials.  Thus, Euler
writes that
\begin{quote}
if~$dy = pdx$ and~$dp = q dx$, then the second differential~$d^2y = q
dx^2$, and so it is clear, as we indicated before, that the second
differential of y has a finite ratio to~$dx^2$.  (ibid., p.~68)
\end{quote}
Here Euler is assuming that~$dx$ involves a \emph{constant}
progression of differentials so that~$ddx=0$.  A more general
situation is dealt with in~\S\,129:
\begin{quote}
129. If the successive values of~$x$, namely,
$x$,~$x^I,x^{II}, x^{III}$,
$x^{IV}, \ldots$, do not form an arithmetic progression, but follow
some other rule, then their first differentials, namely,
$dx$,~$dx^I, dx^{II},\ldots$, will not be equal to each other, and so
we do not have~\hbox{$d^2x=0$}. (ibid., p.~68)
\end{quote}
Such insights are difficult to relate to from the modern viewpoint
centering on the concept of derivative (rather than differential).
This different nature of the infinitesimal calculus as practiced by
both Leibniz and Euler was emphasized in the seminal study \cite{Bos}
(on Leibniz see further in \cite{KS1}, \cite{SK}).  Euler's \S\,138 on
page~72 may therefore come as a surprise to someone trained in
the~$f\rightsquigarrow f'$ tradition:
\begin{quote}
From the word \emph{differential}, which denotes an infinitely small
difference, we derive other names that have come into common
usage. Thus we have the word \emph{differentiate}, which means
\emph{to find a differential}. (ibid., p.~72)
\end{quote}
That's not what the term \emph{differentiate} means today!  But the
real bombshell occurs in Chapter 5, \S\,164.  Euler writes:
\begin{quote}
Let~$p/q$ be a given function whose differential we need to find.
When we substitute~$x + dx$ for~$x$ the quotient becomes
\[
\frac{p+dp}{q+dq}=(p+dp)\left(\frac{1}{q}-\frac{dq}{q^2}\right)=
\frac{p}{q}-\frac{p\,dq}{q^2}+\frac{dp}{q}-\frac{dp\,dq}{q^2}.
\]
When~$p/q$ is subtracted, the differential remains,
$d\left(\frac{p}{q}\right)= \frac{dp}{q}-\frac{p\,dq}{q^2}$,
\end{quote}
and adds:
\begin{quote}
since the term~$\frac{dp\,dq}{q^2}$ vanishes [\emph{ob evanescentem
terminum}~$\frac{dp\,dq}{qq}$ in the original Latin].
\end{quote}
Euler thus obtains the formula for the differential~$d(p/q)$ in terms
of~$dp$ and~$dq$, namely what we would call today the quotient rule.
Note the absence of power series and the presence of the idea of
\emph{discarding higher order terms}, such as~$\frac{dp\,dq}{q^2}$.
We will see that this aspect of Euler's work is at tension with
Edwards' interpretation.  For more details on Euler's foundational
stance, see \cite{Re13}, \cite{Ba14}.  For related updates on Fermat
and Cauchy see respectively \cite{KSS13} and \cite{BK}.

\section{How many strings does a flute have?}

Harold M. Edwards opens his piece on Euler in the \emph{Bulletin AMS}
with a musical metaphor involving the lute and the oud.  These are
similar instruments, but in practice are used to play very different
tunes.  Edwards' metaphor is meant to illustrate a claim concerning
what he refers to as ``Euler's definition of the derivative.''
Edwards claims that, while the tune of Euler's definition may sound
like the ratio of a pair of infinitesimals, in reality something else
is going on:
\begin{quote}
When I understood enough of the context to realize what Euler was
saying, I experienced a \emph{shock} of recognition. It was
practically the same as the definition of the derivative that I
finally chose after decades of teaching calculus: `Rewrite
$\frac{\Delta y}{\Delta x}$ in a way that still makes sense when
$\Delta x = 0$.'%
\footnote{When Edwards proposes to ``rewrite'' the ratio, he does not
mean to replace it by a different quotient but rather replace it by a
different expression which is no longer a quotient.  For example,
$\frac{(x+dx)^2-x^2}{dx} = \frac{(x+dx+x)(x+dx-x)}{dx} = 2x+dx$, and
the last expression \emph{still makes sense} when~$dx=0$.  Now the
phrase ``Rewrite \ldots{} in a way that still makes sense'' is
certainly ambiguous.  Does this entail developing into a power series
and taking the linear term, or perhaps a different technique?  There
seems to be a conflation of \emph{definition of X} and \emph{algorithm
for computing X}, which may be deliberate given Edwards' algorithmic
and constructive philosophical commitments.}
\cite[p.~576]{Ed07b} (emphasis added)
\end{quote}
Alas, we have read the relevant passages in Euler but neither
experienced the epiphanous ``shock'' Edwards reports, nor for that
matter detected any such similarity.

As far as listening to Euler is concerned, at the very least it needs
to be pointed out that Edwards lacks a perfect pitch.  The first false
note is already in his title ``Euler's definition of the derivative,''
for Euler \emph{did not define the derivative} at all.  Therefore the
answer to the question ``What was Euler's definition of the
derivative?'' is: \emph{None}, similarly to the answer to the question
contained in the title of this section.  Namely, the answer is that
Euler doesn't \emph{give} a definition of derivative.

Euler works with differentials throughout (see Section~\ref{one} for
examples).  The differential quotient plays an auxiliary role and
always appears in a relation between differentials, such as the factor
of~$2x$ in~$dy=2x dx$ when~$y=x^2$.  Derivatives don't appear either
in the \emph{Introductio} or in the \emph{Institutiones}, either under
their modern name or as fluxions.  Euler mentions fluxions in \S~115
of the \emph{Institutiones} as the English equivalent of
differentials:
\begin{quote}
The English mathematicians \ldots{} call infinitely small differences,
which we call differentials, \emph{fluxions} and sometimes
\emph{increments}. \cite[\S 115]{Eu55b}
\end{quote}
%
%
The noun \emph{derivative} used in the translation
\cite[\S~235--238]{Eu55b} is Blanton's and doesn't appear in the Latin
original.

Euler's student Lagrange did introduce \emph{la fonction d\'eriv\'ee}
in his article \cite{La72}, but this was well after the publication of
the \emph{Institutiones} \cite{Eu55}, the text Edwards claims as his
source.  One might have thought that Edwards merely simplified the
title for greater accessibility, but the same jarring note is sounded
in the abstract: ``Euler's method of defining the derivative of a
function is not a failed effort to describe a limit.''  Edwards is
still out of tune in his introduction: ``[Euler's] definition of the
derivative is misunderstood primarily because his notion of `function'
is misunderstood.''  \cite[p.~576]{Ed07b}.  A crescendo ``\emph{Of
course} Euler understood limits.  Euler was Euler.  But he rejected
limits as the way to define derivatives'' (ibid.)  is followed by a
coda ``Since the definition of the derivative is still two volumes in
the future,'' (ibid., p.~579).

\section{Where did the infinitesimals go?}

One of Edwards' major mathematical contributions is what is known as
the Edwards curve~$x^2+y^2=1+d\,x^2y^2$, a way of writing certain
elliptic curves that, while apparently less elegant than the
Weierstrass form~$y^2=x^3+ax+b$, turns out to be more efficient,
computationally and constructively speaking.  His text on elliptic
curves \cite{Ed07a} appeared in the same journal a few months earlier.
The link to the Euler text is that, as noted in Edwards' \emph{Essays
in constructive mathematics},
\begin{quote}
Euler too dealt with the curve $y^2= 1 - x^4$ \ldots, for which
explicit and beautiful formulas can be developed for the addition law,
\ldots{} To require that it be put in Weierstrass normal form before
the group law is described loses certain symmetries that deserve to be
kept. \cite[p.~127]{Ed05}.
\end{quote}
This is a fascinating historical observation, but reading Edwards on
Eulerian calculus one may not even suspect that the~$\omega$ appearing
in Euler's discussion of the differential ratio is infinitesimal.  In
fact, \cite{Ed07b} is nearly an infinitesimal-free zone (though he
does mention an infinitely small~$z$ on page 578).  His comment on
students being ``taught to shrink from differentials as from an
infectious disease'' \cite[p.~579, note~5]{Ed07b} appears to share
epidemiological concerns with Cantor's vilification of infinitesimals
as the \emph{cholera bacillus of mathematics} \cite[p.~505]{Me}.
However, none of this is faithful to Euler, as we saw in
Section~\ref{one}.  In fact Edwards' comment is but the tip of the
iceberg of attempted constructivist deconstructions of infinitesimals;
see \cite{Ka15} for further details.

Now it is certainly possible, mathematically speaking, to redefine the
derivative as the coefficient of the linear term in the Taylor series.
However, this was Lagrange's approach, not Euler's; see for example
\cite[p.~100]{Gr}.  Edwards, possibly due to his constructivist
leanings, appears to favor Lagrange's definition over Euler's.
Edwards apparently believes that any great mathematician would agree
with him on this point, and moreover he believes that Euler was a
great mathematician (``Euler was Euler'').  Edwards' syllogism,
however, fails to persuade someone who wishes to undestand Euler's
actual procedures, rather than what they should ``really mean'' for
someone who thinks he knows what the unique consistent interpretation
of this must be so as to save Euler's honor.  Readers of Edwards'
second \emph{Bulletin} article presumably expected to get a fair
picture of Euler's foundational stance from reading the article.  This
they arguably did \emph{not} get.

Edwards' article might have been more appropriately titled
\emph{Lagrange's definition of the derivative}.  Note however that in
the second edition of his \emph{M\'ecanique Analytique}, Lagrange
fully embraced infinitesimals in the following terms:
\begin{quote}
Once one has duly captured the spirit of this system [i.e.,
infinitesimal calculus], and has convinced oneself of the correctness
of its results by means of the geometric method of the prime and
ultimate ratios, or by means of the analytic method of derivatives,
one can then exploit the infinitely small as a reliable and convenient
tool so as to shorten and simplify proofs'' \cite[p.~iv]{Lag}.
(translation ours)
\end{quote}
While Lagrange didn't make the shortlist in \cite[p.~105]{Ed81}, we
highly recommend both Euler and Lagrange.  \emph{Read the masters!}

\subsection{Postscript}

In a response to this article, H. Edwards finds our ideas strange,
claims that they derive from set theory, and furthermore that the
ideas of nonstandard analysis are strange; see \cite{Ed15}.  However,
our article does not say \emph{a single word} about either set theory
or nonstandard analysis.  Edwards brings up the subject with
apparently little justification.  What purpose is served by claiming
that nonstandard analysis is strange?  Apparently the only purpose is
an attempt to divert attention from the shortcomings of his own
reading of Euler through a demagogical ad hominem attack, unworthy of
a response.  As a postscript we will merely point out that Terence Tao
has recently published a number of monographs where he exploits
ultraproducts in general, and Robinson's infinitesimals in particular,
as a fundamental tool; see e.g., \cite{Ta14}.

\section*{Acknowledgments}

The work of Vladimir Kanovei was partially supported by the Russian
Scientific Fund (project no. 14-50-00150) and by RFBR grant
13-01-00006.  M.~Katz was partially funded by the Israel Science
Foundation grant no.~1517/12.  The influence of Hilton Kramer
(1928-2012) is obvious.

\medskip\noindent \textbf{Vladimir Kanovei} graduated in 1973 from
Moscow State University.  He obtained a Ph.D. in physics and
mathematics from Moscow State University in 1976.  In 1986, he became
Doctor of science in physics and mathematics at Moscow Steklov
Mathematical Institute (MIAN).  He is currently Leading Researcher at
the Institute for Information Transmission Problems (IITP), Moscow,
Russia, 

\noindent
\url{http://www.iitp.ru/en/users/156.htm} He currently holds a
part-time job as Professor at Moscow State University of Railway
Engineering (MIIT), Moscow, Russia.  His research interests in
mathematics include set theory, logic, nonstandard analysis, history
of infinitesimals.  His selected monographs are:

1) (with M.Reeken) Nonstandard analysis, axiomatically. Springer
Monographs in Mathematics 2004, ISBN: 978-3-540-22243-9

2) Borel equivalence relations, structure and
classification. University Lecture Series 44. American Mathematical
Society, Providence, RI, 2008, ISBN: 978-0-8218-4453-3

3) (with M.Sabok and J.Zapletal) Canonical Ramsey Theory on Polish
Spaces. Cambridge University Press, Cambridge, 2013, ISBN
978-1-107-02685-8

\medskip\noindent{\bf Karin Usadi Katz} (B.A. Bryn Mawr College, '82;
Ph.D. Indiana University, '91) teaches mathematics at Bar Ilan
University, Ramat Gan, Israel.

\medskip\noindent \textbf{Mikhail G. Katz} (B.A. Harvard University,
'80; Ph.D. Columbia University, '84) is Professor of Mathematics at
Bar Ilan University.  In the fall of the '14-'15 academic year, he
taught, together with a colleague, True Infinitesimal Calculus using
H. J. Keisler's book \emph{Elementary Calculus} to a group of 120
freshmen, prompting one of the students to comment as follows in a
somewhat Hebraized English: ``I really enjoyed learning this course
and it was much understandable and `fun' from the standard approach.
My friends in other universities didn't really enjoy Calculus and
actually kind of hated it but we really enjoyed learning it and even
had some jokes about infinitesimal things in our lives :)''

\medskip\noindent \textbf{David Sherry} is fortunate to be Professor
of Philosophy at Northern Arizona University, Flagstaff, AZ, in the
cool pines of the Colorado Plateau.  E-mail: david.sherry@nau.edu.  He
took his Ph.D. at Claremont Graduate School in 1982.  His research
interests include history and philosophy of mathematics, especially
applied mathematics, and history of modern philosophy (16th-18th
century).  He regularly offers courses in history of modern
philosophy, symbolic logic, and metaphysics.  When his nose is not in
a book, he's likely hiking and looking at birds with his wife, making
music, gardening, or playing with his grandchildren.  Fortunate
indeed.

\end{document}